\documentclass[a4paper,12pt,leqno]{article}

\usepackage{graphics}
\usepackage[T1]{fontenc}
\usepackage{lineno}
\usepackage[utf8]{inputenc}
\usepackage{amsmath,mathrsfs,amssymb,amsfonts,nicefrac}
\usepackage[english]{babel}
\usepackage[active]{srcltx}
\usepackage{tikz}
\usepackage{dsfont}
\usetikzlibrary{patterns}

\numberwithin{equation}{section}

\newtheorem{thm}{Theorem}[section]

\newtheorem{lm}[thm]{Lemma}  

\newtheorem{rem}[thm]{Remark}  

\newcommand{\red}{\textcolor{red}}
\newcommand{\RR}{\mathbb{R}}   
\newcommand{\R}{\mathbb{R}}   
\newcommand{\di}{\displaystyle}

\renewcommand{\epsilon}{\varepsilon}
\newcommand{\e}{\varepsilon}
\def\un{{\mathbf{1}}}
\def\o{\overline}

\begin{document}   
   
\title{\textbf{Improved bounds for reaction-diffusion propagation  with a line of nonlocal diffusion}}

\author{
{\bf Anne-Charline Chalmin} \\ 
Institut de Math\'ematiques de Toulouse (UMR CNRS 5219) \\ 
Universit\'e Toulouse III,
118 route de Narbonne\\
31062 Toulouse cedex, France \\ 
\texttt{abchalmin@gmail.com}\\ 
\\[2mm]
{\bf Jean-Michel Roquejoffre} \\ 
Institut de Math\'ematiques de Toulouse (UMR CNRS 5219) \\ 
Universit\'e Toulouse III,
118 route de Narbonne\\
31062 Toulouse cedex, France \\
\texttt{jean-michel.roquejoffre@math.univ-toulouse.fr}}
   
\maketitle   
\begin{abstract} We consider here a model of accelerating fronts, introduced in \cite{BCRR}, consisting of one equation with nonlocal diffusion on a line, coupled via the boundary condition with a reaction-diffusion equation of the Fisher-KPP type in the upper half-plane.  It was proved in \cite{BCRR} that the propagation is accelerated in the direction of the line exponentially fast in time. We make  this estimate more precise by computing an explicit correction that is algebraic in time. Unexpectedly, the solution mimicks the behaviour of the solution of the  equation linearised around the rest state 0 in a closer way than in the classical fractional Fisher-KPP model.
\end{abstract}

\begin{center}
{\it This paper is dedicated to S. Salsa, as the expression of our friendship and respect.}
\end{center}

\section{Introduction}    
\subsection{Model and question}
Consider the following system, with unknowns $u(t,x)$ and $v(t,x,y)$, where $(t,x,y)\in\RR_+\times\RR\times\RR_+$
\begin{equation}
\label{e1.1}
\left\{\begin{array}{rcll}
\partial_t v-\Delta v &=&f(v),&\quad t>0,x \in\R,\ y >0, \\
\partial_tu+(-\partial _{xx})^{\alpha}u&=& -\mu u +\nu v&\quad t>0, x \in \R,\ y=0,\ 
t>0\\
-\partial_y v&=&\mu u - \nu v, &\quad t>0, x \in \R,\ y=0.\\
\end{array}
\right.
\end{equation}
The real number $\mu$ is a positive given parameter, and the nonlinear term $f$ is chosen as
$$
f(v)=av-g(v),
$$
with $a>0$ and $g$ of class $C^2$, with $g\geq0$, convex, $g(0)=g'(0)=0$, $g'(+\infty)>a$. The equation for $v$ in the upper half-plane is therefore a variant of the Fisher-KPP equation, in reference to the 
pioneering works of Fisher \cite{F} and of Kolmogorov, Petrovskii and Piskunov \cite{KPP}.
 The operator $(-\partial_{xx})^\alpha$ is the fractional Laplacian of order $\alpha\in(0,1)$:
 $$
 (-\partial_{xx})^\alpha u(x)=c_\alpha\ {\mathrm{P.V.}}\biggl(\int_{\RR}\frac{u(y)-u(x)}{\vert x-y\vert^{1+2\alpha}}dy\biggl),
  $$
the constant $c_\alpha>0$ being chosen so that the symbol of $(-\partial_{xx})^\alpha$ is $\vert\xi\vert^{1+2\alpha}$.
The initial datum is chosen as
\begin{equation}
\label{e1.2}
(u(0,x),v(0,x,y))=(\delta_0\un_{(-x_0,x_0)}(x),0)
\end{equation}
where $x_0$ and $\delta_0$ are given positive constants. Their value is not relevant for the discussion, one may think them as small. Under the listed assumptions, system \eqref{e1.1}  has a unique global smooth solution, that is also globally bounded as well as its derivatives, see \cite{BCRR}. The question under study is the behaviour of $(u(t,x),v(t,x,y))$ for large $t$.
\subsection{Motivation, context, known results}
System \eqref{e1.1} is relevant in the study of the influence of a line having a fast diffusion of its own, that exchanges with an adjacent domain of the plane (here, the upper half plane), in which reactive and diffusive phenomena occur. The application is the modelling of how biological invasions can be enhanced by transportation networks, see \cite{BCRR_Schwartz} for an overview. In this context, $u(t,x)$ represents the density of individuals on the line, and $v(t,x,y)$ represents the density of individuals in the upper half-plane. Exchanges occur through the Robin condition $-\partial_yv(t,x,0)=\mu u(t,x)-v(t,x,0)$. 

System \eqref{e1.1} was first introduced by H. Berestycki, L. Rossi and the second author in \cite{BRR2}. There, the diffusion on the line (that we called "the road", while the upper half plane was called "the field") took the form $-D\partial_{xx}$, with $D>0$, possibly large. The effect of the line may be accounted for as follows: when not present, the model amounts to the single Fisher-KPP equation with unknown $v(t,X)$, $X\in\RR^2$:
\begin{equation}
\label{e1.3}
\left\{\begin{array}{rll}
v_t-\Delta v=&f(v),\quad t>0,\ X\in\RR^2\\
v(0,X)=&\delta_0\un_{(-x_0,x_0)^2}(X), \quad X\in \RR^2.
\end{array}\right.
\end{equation}
Note that here, we need to shift the mass from the line to the plane in order to avoid the trivial solution $v\equiv0$. We have (Aronson, Weinberger \cite{AW})
\begin{equation}
\label{e1.4}
\begin{array}{rll}
&\hbox{for all $\e>0$,}\quad \di\lim_{t\to+\infty}\inf_{\vert X\vert\leq(c_*-\e)t}v(t,X)=v_0\\
&\hbox{for all $\e>0$,}\quad \di\lim_{t\to+\infty}\sup_{\vert X\vert\geq(c_*+\e)t}v(t,X)=0,\\
\end{array}
\end{equation}
where $c_*=2\sqrt{a}$, and $v_0$ is the unique positive zero of $f$, whose existence is granted by the assumptions. In other words, the stable state $v_0$ invades the whole space at speed $c_*$.
Reverting to \eqref{e1.1}, and concentrating on what happens on the line (or its vicinity), first when the diffusion is $-D\partial_{xx}$, then when it is $(-\partial_{xx})^\alpha$. In the first case, the main result of \cite{BRR2} is the existence of $c_*(D)>0$, with $\di\liminf\frac{c_*(D)}{\sqrt D}>0$, such that invasion occurs at speed $c_*(D)$ on the line and in the upper half plane, at finite distance from the line. This ashows the importance of the line on the overall propagation. The limiting states for $u$ and $v$ are $u_\infty\equiv\di\frac{v_0}\mu$, $v_\infty\equiv v_0$, a property that is not entirely trivial, and also proved in \cite{BRR2}. 

 The effect of the nonlocal diffusion $(-\partial_{xx})^\alpha$ was studied for the first time in \cite{BCRR} by H. Berestycki, L. Rossi and the two authors of the present paper. The main result of \cite{BCRR} is the following.
\begin{thm}
\label{t1.0}
Define $\lambda_*=\di\frac{a}{1+2\alpha}.$ Then we have
\begin{equation}
\label{e1.5}
\begin{array}{rll}
&\hbox{for all $\e>0$,}\quad \di\lim_{t\to+\infty}\inf_{\vert x\vert\leq e^{(\lambda_*-\e)t}}(u(t,x),v(t,x,y))=\left(\di\frac{v_0}\mu,v_0\right)\\
&\hbox{for all $\e>0$,}\quad \di\lim_{t\to+\infty}\sup_{\vert x\vert\geq e^{(\lambda_*+\e)t}}(u(t,x),v(t,x,y))=0.\\
\end{array}
\end{equation}
In \eqref{e1.5}, the limits of $v$ should be understood pointwise in $y$. 
\end{thm}
Let us note that this result may be parallelled by the following one: let us bluntly replace the exchange term $\mu u-v$ in the equation for $u$ by the reaction term $f(u)$ (so that we shift the whole weight of the reaction from the upper half plane to the line), so as to obtain
\begin{equation}
\label{e1.6}
\left\{\begin{array}{rll}
u_t+(-\partial_{xx})^\alpha u=&f(u)\quad (t>0,\ x\in\RR)\\
u(0,x)=&\delta_0\un_{(-x_0,x_0)}(x).
\end{array}\right.
\end{equation}
Then, Cabr\'e and the second author \cite{CR} proved that invasion at the same rate as in Theorem \ref{t1.0} occurs. Thus, $u(t,x)$ actually behaves just as in equation \eqref{e1.6} at the leading order.

While Theorem \ref{t1.0} captures the essence of the main features of the invasion phenomenon, it is interesting to ask whether the asymptotics can be made a little more precise. Indeed there is, in Theorem \ref{t1.0}, a lot of room between the upper and lower bound. For instance the level sets of $u$ may advance like $t^pe^{\lambda_* t}$, where $p$ could be any real number.  
This question can also be asked for the simpler model \eqref{e1.6}, all the more as one may give the following heuristics : the dynamics of \eqref{e1.6} being driven by the small values of $u$ (given the concavity of $u$ they are, loosely speaking, the most unstable ones in the range of $f$), so that the dynamics of the level sets is really given by the linear equation
$$
u_t+(-\partial_{xx})^\alpha u=au.
$$
Call $G_\alpha(t,x)$ the fractional heat kernel, we have $G_\alpha(t,x)\lesssim\di\frac{t}{\vert x\vert^{1+2\alpha}}$ for large $t$ and $x$, see \cite{Kolo} for instance. Then we have
$$
u(t,x)\lesssim \frac{te^{at}}{\vert x\vert^{1+2\alpha}}, 
$$
still  for large $t$ and $x$. So,  a level set of $u$ will move like $t^{\frac{1}{1+2\alpha}}e^{\frac{at}{1+2\alpha}}$. This heuristics does not give the correct sharper behaviour, as was proved by Cabr\'e and the two authors of the paper \cite{CCR}: a level set $\{x(t)\}$ of $u$ will in fact be such that 
$\vert x(t)\vert e^{-\frac{at}{1+2\alpha}}$ is bounded, that is, there is no polynomial correction in the expansion of $x(t)$. 

Consider now the  linearised version of \eqref{e1.1}:
\begin{equation}
\label{e1.7}
\left\{\begin{array}{rcll}
\partial_t v-\Delta v &=&av,&\quad t>0,x \in\R,\ y >0, \\
\partial_tu+(-\partial _{xx})^{\alpha}u&=& -\mu u +\nu v&\quad t>0, x \in \R,\ y=0,\ 
t>0\\
-\partial_y v&=&\mu u - \nu v, &\quad t>0, x \in \R,\ y=0.\\
\end{array}
\right.
\end{equation}
Let us call this time $G_\alpha(t,x)$ the solution $u(t,x)$ with the initial datum $\delta_{x=0}$, that is, the $u$-component of the fundamental solution. Then the first author proved \cite{ACC} (a more precise estimate will be stated later).
$$
G_\alpha(t,x)\lesssim \frac{e^{at}}{t^{\frac32}\vert x\vert^{1+2\alpha}} ,\quad t\to+\infty,\ \vert x\vert\to+\infty.
$$
And so, a level set $\{x(t)\}$ of the solution $u(t,x)$ of \eqref{e1.7} will move like $t^{-\frac{3}{2(1+2\alpha)}}e^{\frac{at}{1+2\alpha}}$. The question that we want to address in this paper is whether a discrepancy of the same kind holds between the linear and nonlinear equation.
\subsection{Result and organisation of the paper}
Surprisingly, and in contrast to what happens with \eqref{e1.6}, the linear equation \eqref{e1.7} mimicks the behaviour of the nonlinear one \eqref{e1.1} in a better fashion than for the fractional Fisher-KPP equation. The result that we are going to prove is the following.
\begin{thm}
\label{t1.1} 
Consider any $\lambda\in\left(0,\di\frac{v_0}\mu\right)$. Let $x_\lambda(t)$ be the largest $x$ such that $u(t,x)=\lambda$ or $u(t,-x)=\lambda$. Then, for all $\delta>0$, there is $T_{\lambda,\delta}>0$ such that, for all $t\geq T_{\lambda,\delta}$ we have 
\begin{equation}
\label{e1.8}
\frac{e^{\frac{a}{1+2\alpha}t}}{t^{\frac3{2(1+2\alpha)}+\delta}}\leq x(t)\leq\frac{e^{\frac{a}{1+2\alpha}t}}{t^{\frac3{2(1+2\alpha)}-\delta}}
\end{equation}
\end{thm}
In fact, the upper bound is more precise, as we may choose $\delta=0$ there. To improve the lower bound seems to us more challenging, and will be addressed in a future work.

The paper is organised as follows. In Section \ref{s2}, we explain the strategy of the proof of Theorem \ref{t1.1} and discuss some perspectives that our work has opened. In Section \ref{s3} we address the underlying mechanism of Theorem \ref{t1.1}, namely, the transients of the one-dimensional Fisher-KPP equation with Dirichlet boundary conditions, this is a result of independent interest for the Fisher-KPP equation. We then devote a short section to quantify how the exchanges between the road and the field are organised. The proof of Theorem \ref{t1.1} is then displayed in Section 5. In the whole paper, the computations will be greatly simplified when we take a function  $g\geq0$, smooth, convex, supported in 
$(\theta,1]$ for some $\theta\in(0,1)$, with $g(1)=1$. Threfore, the computations will sometimes be carried out with this type of nonlinearity in order to highlight the main ideas,  before being extended  to the general Fisher-KPP nonlinearity. Also, from now on we will assume, without loss of generality, that $a=1$.

\section{The underlying mechanism of Theorem \ref{t1.1}, discussion}\label{s2}
The starting point of this paper was the following numerical simulations, carried out in the PhD thesis of the first author \cite{ACC}.
\begin{figure}[h]
	\includegraphics[width=0.32\textwidth]{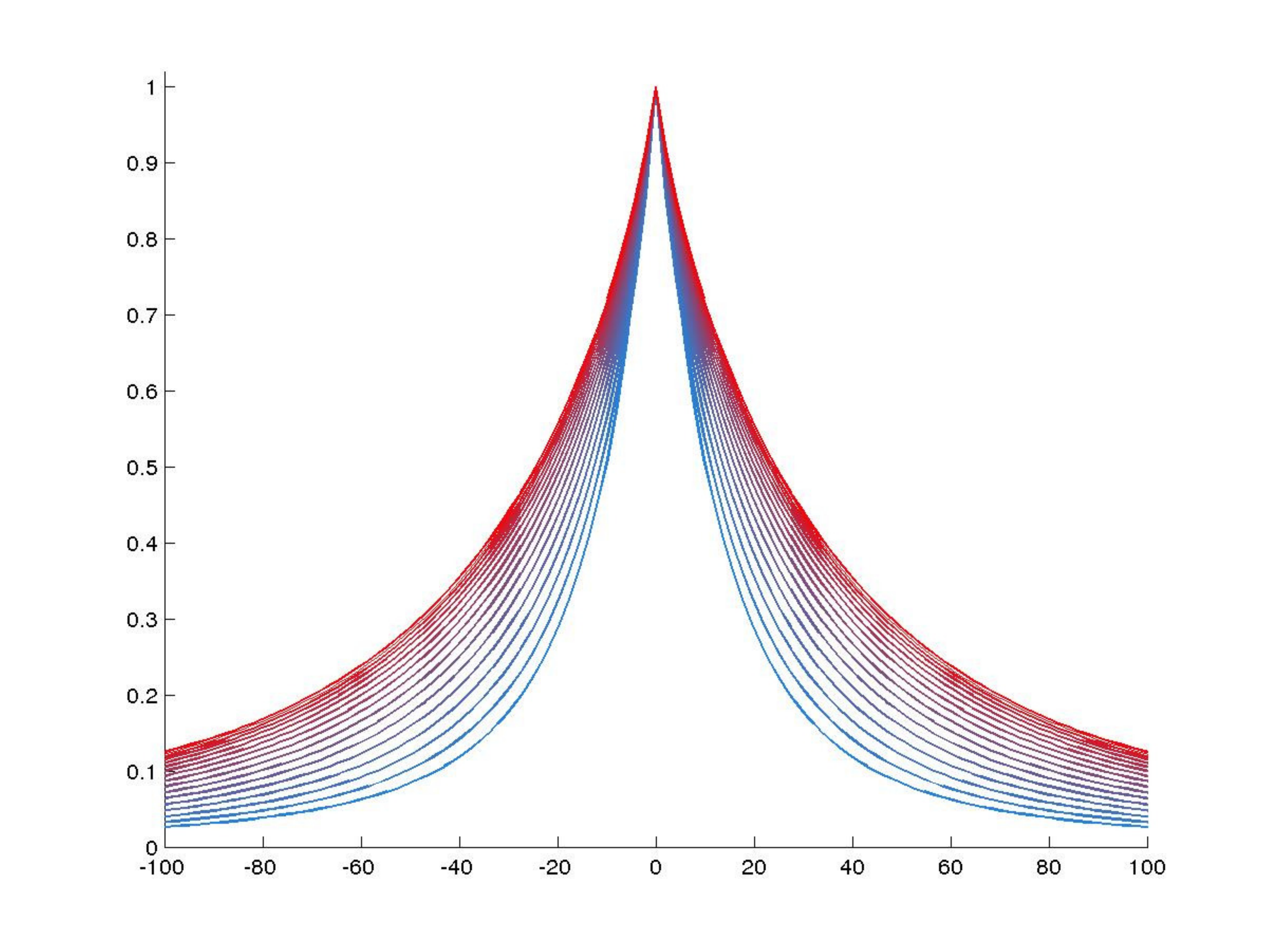}
	\includegraphics[width=0.32\textwidth]{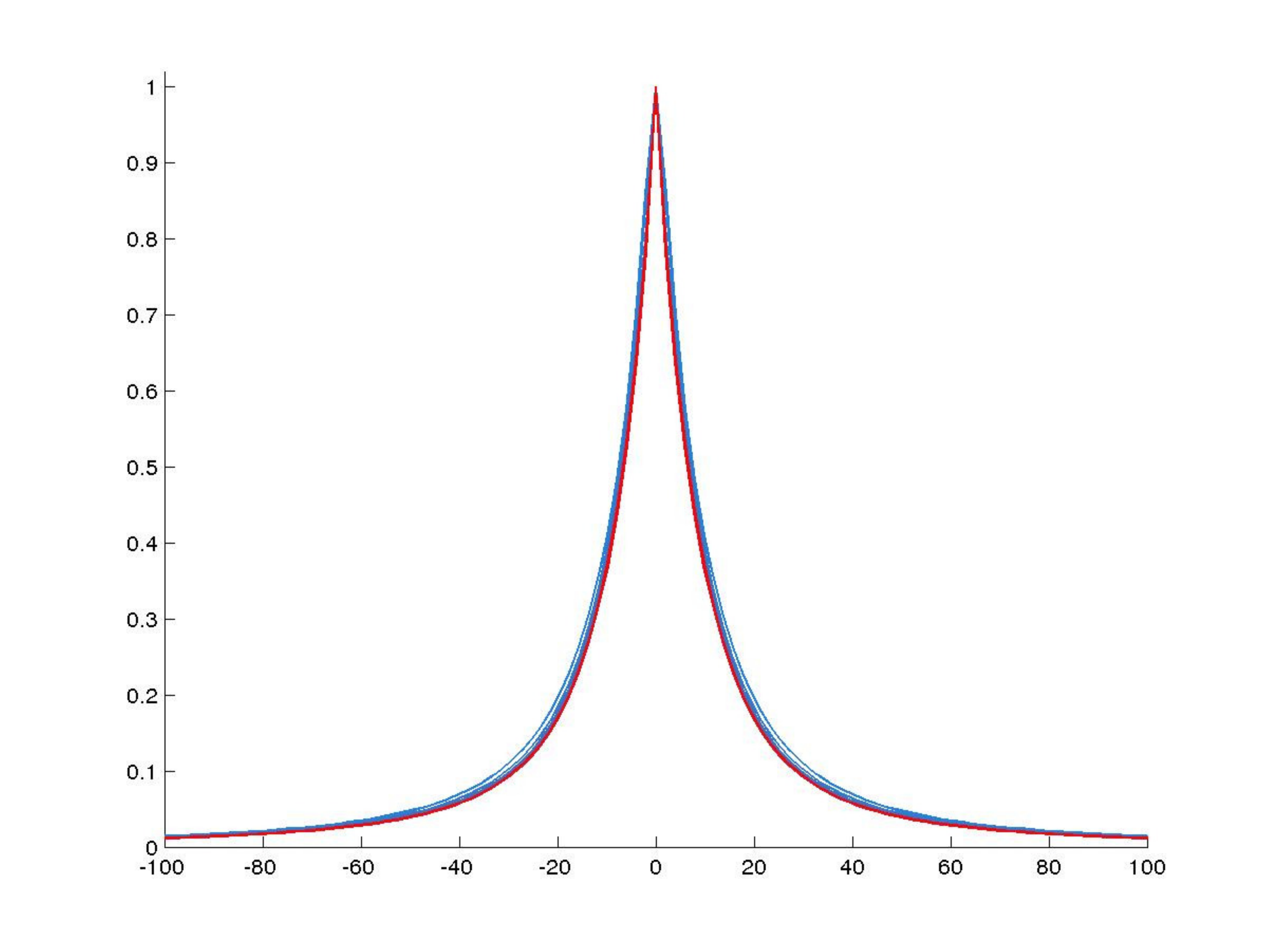}
\includegraphics[width=0.32\textwidth]{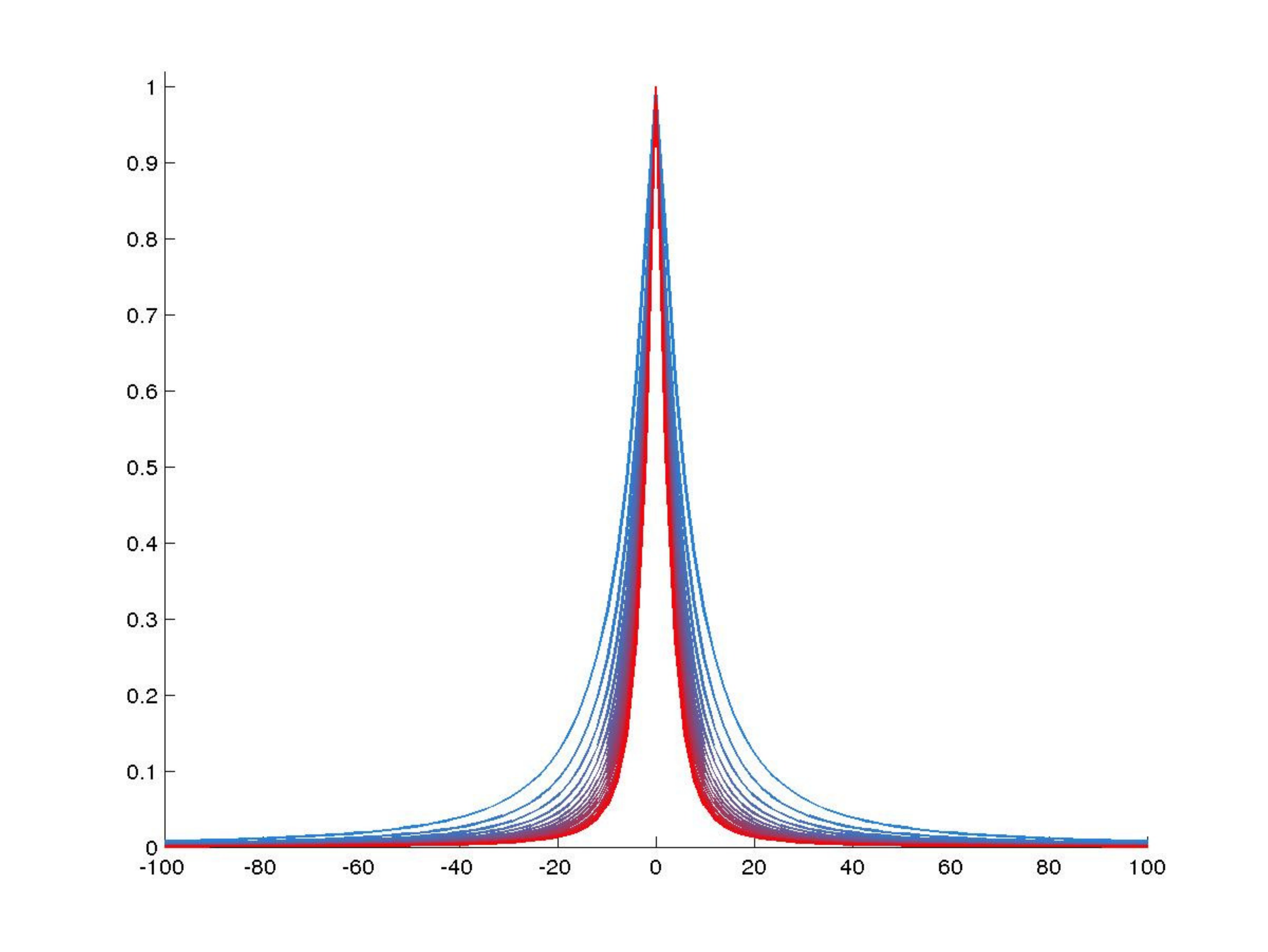}
\caption{The different renormalisations}
\end{figure}
The figure represents the evolution of $u(t,t^{-m}e^{\frac{t}{1+2\alpha}})$, with, from left to right: 
$m=0$, $m=\di\frac{3}{2(1+2\alpha)}$ and $m=\di\frac{3}{1+2\alpha}$. The gradation of colours from blue to red represents the advance in time, blue standing for the earlier stages of the development. One clearly sees the stabilisation mechanism for the middle value of $m$, and this came to us as a surprise.  However, this suggests the following idea: the $t^{-3/2}$ term being typical of the one-dimensional Dirichlet heat equation, we thought that it was interesting to understand this feature in a little more depth.

 Assuming - which will turn out to be  a good approximation - that $\partial_{xx}v$ is small, this suggests in fact that, for a fixed $x$, the function $v(t,x,y)$ behaves like a solution of the one dimensional Fisher-KPP equation
 $$
 \begin{array}{rll}
 w_t-w_{yy}=&f(v),\quad t>0,\ y>0\\
 w(t,0)=&0,
 \end{array}
 $$
This is even more evident when one takes $f(v)=av-g(v)$, $g$ vanishing on a small interval to the right of 0. The initial value (or, at least, the value of $v$ at any small positive time) is small, dictated by the size of $u(1,x)$.  The Dirichlet boundary condition is the most convenient one that allows to put below the solution $v(t,x,y)$ of \eqref{e1.1} a barrier devised on the model of $w(t,y)$, with an initial datum suitably dictated by the behaviour of $u(t,x)$ (the solution of \eqref{e1.1} on the line) at infinity. Of course, with this particular condition, the role of the line seems to be forgotten, such is not exactly the case, as long as we prove - as will be done in the course of this work - some easy lemmas that describe how the communication between the road and the field is organised. 

Let us briefly discuss the optimality of our estimates. Of course, the corrections of the exponents in \eqref{e1.8} by a small $\delta>0$ shows that there is still a room for improvement.  In particular, one could ask whether replacing the Dirichlet boundary condition by the exchange condition $-\partial_yv+\nu v$ in the 1D Fisher-KPP equation would lead to the optimal bounds. In fact, the best strategy would probably be to investigate the full one-dimensional problem with unknowns $(u(t),v(t,y))$
$$
\left\{
\begin{array}{rll}
v_t-v_{yy}=&f(v),\quad t>0,\ y>0\\
-v_y(t,0)=&\mu u(t)-\nu v(t,0)\quad t>0\\
\dot u(t)=&\nu v(t,0)-\mu u(t).
\end{array}
\right.
$$
We choose no to do it here, as it would involve, in our opinion, heavier computations with possibly no real further understanding of the mechanisms at work. So, we leave this task for a future contribution.
\section{The transients for 1D Fisher-KPP propagation with small initial data}\label{s3}
In this section we consider a function $g\geq0$, smooth, convex, supported in 
$(\theta,1]$, with $g(1)=1$ (this last assumption is in fact unnecessary). Pick a small $\e>0$ {with $\e < \theta$}. The goal of this section is to understand how much time
it will take to the solution of the model Fisher-KPP equation
\begin{equation}
\label{e2.1}
\left\{\begin{array}{rll}
v_t-v_{yy}=& v-g(v)\quad(t>0,\ y\geq0)\\
v(0,y)=&\e\un_{[1/2,1]}(y)\\
v(t,0)=&0
\end{array}\right.
\end{equation}
to reach the value $\theta$ at finite distance from $y=0$. First, let us note that the value $\theta$ will eventually be reached at a finite distance (both with respect to $t$ and $\e$) from $y=0$. Indeed, a classical sub-solution argument (see for instance Berestycki-Hamel-Roques \cite{BHR}) implies that $v$ will converge to the unique nontrivial solution $v_\infty$ of 
\begin{equation}
\label{e2.2}
\begin{array}{rll}
-v_\infty''=&v_\infty'-g(v_\infty')\quad(y>0)\\
v_\infty'(0)=&0,
\end{array}
\end{equation}
 which satisfies $v_\infty'(0)>0$, hence uniformly bounded from below on every set of the form $[y_0,+\infty)$, $y_0>0$.
On the other hand, as
$\e\to0$, the time that it will take to $v$ to come close to $v_\infty$ will grow infinitely, and our aim is to devise an upper bound that will be precise up to algebraic powers of $\e$.
\begin{thm}
\label{t2.1}
Let $v_\e$ be the solution to \eqref{e2.1}, and ${\e<}\lambda<v_\infty(1)$. Define $T_\e$ as the first time $t$ such that 
\begin{equation}
\label{e2.3}
v_\e(t,1)=\lambda.
\end{equation}
Then, for all $\delta>0$, there is $Q_\delta>0$, possibly blowing up as $\delta\to0$, such that 
\begin{equation}
\label{e2.4}
\frac{e^{T_\e}}{T_\e^{\frac32+\delta}}\leq\frac{Q_\delta}\e.
\end{equation}
\end{thm}
It is worth saying a word on the scenario leading to \eqref{e2.4}, and the special structure of the nonlinearity
$f(v)={v}-g(v)$ will make it especially obvious: the region where the solution will first reach a nontrivial value is 
not close to 0, but at a large distance from 0. At this stage, one could think of invoking classical results on Fisher-KPP
propagation for studying how much more time $v_\e$ will take to be nontrivial near $y=0$. This is not the correct intuition, 
because it would lead to a $T_\e$ that would be largely over-estimated. The mechanism is in fact closer to that of nonlocal 
Fisher-KPP propagation \cite{CR}, \cite{CCR}. It is also not so far from what happens with the classical Fisher-KPP with 
slowly decreasing initial data, \cite{H}, \cite{HR}.

\noindent{\bf Proof of Theorem \ref{t2.1}.} For small, or, even, finite $t$, (for instance $t\in[1,2]$) we have $v(t,y)<\theta$ as soon as $\e>0$ is small enough. Let us 
make this assumption; as soon as $v_\e\leq\theta$ everywhere we have $g(v_\e)\equiv0$ and, thus:
\begin{equation}
\label{e2.6}
v_\e(t,y)=\e e^{t}\int_{\frac 12}^1\frac{e^{-\frac{(y-y')^2}{4t} }-e^{-\frac{(y+y')^2}{4t}}}{\sqrt{4\pi t}}dy'.
\end{equation}
For $t\geq1$ and $y'\in[ {\frac12},1]$ we have for all $y>1$
$$
 e^{-\frac{(y-y')^2}{4t} }-e^{-\frac{(y+y')^2}{4t}}\leq C\frac{yy'}te^{-\frac{y^2}{5t}},
$$
$C>0$ universal. So we have, for a possibly different $C>0$:
\begin{equation}
\label{e2.5}
v_\e(t,y)\leq C\e \frac{e^{t}}t.\frac{y}{\sqrt t}e^{-\frac{y^2}{5t}}.
\end{equation}
For 	all fixed $t\geq2$, the maximum in $y$ of the right handside of \eqref{e2.5} is 
taken at 
$$
y=z_0\sqrt t,
$$
where $z_0$ is the maximum of $m(z):=ze^{-z^2/5}$. Call $m_0$ the (easily computable) maximum of
$m$, a sufficient condition to have $v(t,y)\leq\theta$ everywhere is to have, from \eqref{e2.5}:
$$
Cm_0\e \frac{e^{at}}t\leq\theta. 
$$
Define $T_\e^1\geq 2$ as
\begin{equation}
\label{e2.10}
\frac{e^{T_\e^1}}{T_\e^1}=\frac{\theta}{m_0C\e}.
\end{equation}
So, we have easily proved that $v_\e$ reaches a nontrivial value in a time of (roughly) the order of ${\mathrm{ln}}(\e^{-1})$, but this value is reached at 
$y\sim \biggl({\mathrm{ln}}\e^{-1}\biggl)^{1/2}$. To study what happens at finite distance to $y=0$, consider $L>0$ large. There is $c_L>0$, with, in the worst case scenario
$$
\lim_{L\to+\infty}c_L=0,
$$
such that, in the limit $\e\to0$, we have, from \eqref{e2.6}:
\begin{equation}
\label{e2.7}
v_\e(T_\e^1,y)\geq\di\frac{ {\e}c_L}{\sqrt{T_\e^1}},\ \hbox{for $1\leq y\leq L$.}
\end{equation}
Let $e_L(y)$ be the first Dirichlet eigenfunction of $-\partial_{xx}$ on 
$(1,L)$, we normalise so that its maximum is 1. Thus we have
$$
e_L(y)=\sin\biggl(\frac{\pi}{L-1}(y-1)\biggl),
$$
with first eigenvalue
$$
\lambda_1(L)=\frac{\pi^2}{(L-1)^2}.
$$
Let $\underline v_{\e,L}(t,y)$ solve
\begin{equation}
\label{e2.8}
\begin{array}{rll}
(\partial_t-\partial_{xx}-1)\underline v_{\e,L}=&0\quad (t>T_\e^1,\ y\in(1,L))\\
\underline v_{\e,L}(t,1)=\underline v_{\e,L}(t,L)=&0\quad (t\geq T^1_\e)\\
\underline v_{\e,L}( {T_\e^1},y)=&\di\frac{ {\e}c_L}{\sqrt{T_\e^1}}e_L(y).
\end{array}
\end{equation}
On the one hand, we have 
\begin{equation}
\label{e2.9}
\underline v_{\e,L}(t,y)=\frac{{\e}c_Le^{(1-\lambda_1(L))(t-T_\e^1)}}{\sqrt{T_1^\e}}e_L(y).
\end{equation}
On the other hand we have $\underline v_{\e,L}(t,y)\leq v_{\e}(t,y)$ for $t\geq T_\e^1$, as long as 
$\underline v_{\e,L}(t,y)$ is globally less than $\theta$. This last condition is fulfilled as long as
\begin{equation}
\label{e2.11}
t-T_\e^1\leq\frac{1}{2(1-\lambda_1(L))}{\mathrm{ln}}\left(\frac{\theta {^2}T_1^\e}{ {c_L}}\right),
\end{equation}
and the maximum is exactly $\theta$ at equality. Note that it is attained far away from the origin, that is, at 
$y_L=\di\frac{L-1}2 {+1}$. However, the situation is not as bad as before, because we now have
$$
e_L(y)\sim \frac{\pi}{L-1}(y-1),\ \hbox{for $y-1<\!<L$}.
$$
This is certainly a small quantity, but it is independent of $\e$. Let us set 
\begin{equation}
\label{e2.12}
T_\e^2=T_\e^1+\frac{1}{2(1-\lambda_1(L))}{\mathrm{ln}}(\theta^2 T_\e^1).
\end{equation}
We have now, from \eqref{e2.11}:
$$
v_\e(T_\e^2,y)\geq\frac{\theta\pi}{ {L-1}}(y-1), \quad \text{for $y<\!<L$.} 
$$
From now on, once again by a classical sub-solution argument, there is $\tilde T_L>0$ (independent of $\e$), blowing up as $L\to+\infty$, such that
$$
v_\e(T_\e^1+T_\e^2+\tilde T_L,{2})=\lambda.
$$
This is not exactly \eqref{e2.3}, but we are now quite close to it: from the Harnack inequality we have 
$$
v_\e(T_\e^1+T_\e^2+\tilde T_L+1,{1})\geq q\lambda,
$$
for some universal $q>0$, and the same sub-solution argument yields the \eqref{e2.3}, at a time of the form $T_\e^1+T_\e^2+\tilde T_L+\tilde T_L'$, the new constant $\tilde T_L'$ being $\e$-indeendent. Set $T_\e=T_\e^1+T_\e^2+\tilde T_L+\tilde T_L'$; it now suffices to notice that \eqref{e2.11} implies that
$$
T_\e=T_\e^1+\frac12\left(1+O\left(\frac1{L^2}\right)\right){\mathrm{ln}}\left(\frac{\theta^2T_1^\e}{ {c_L}}\right)+\tilde T_L+\tilde T_L',
$$
which, combined to \eqref{e2.10}, implies 
$$
\frac1\e=\frac{ {\sqrt{c_L}}e^{T_\e-\tilde T_L}}{\theta^{1+O(1/L^2)}T_\e^{3/2+O(1/L^2)}(1+T^{-1}_\e{\mathrm{ln}}T_\e-T^{-1}_\e(\tilde T_L+\tilde T_L'))}.
$$
Let us denote by $\mu_L$ a common bound for the two $O(\di\frac1{L^2})$ appearing in the above expression. We now pick a small $\delta$ and choose $L>0$, denoted by $L_\delta$, such that
$$
\mu_{L_\delta}=\delta,\quad Q_\delta=\frac{e^{\tilde T_{L_\delta}}}{\theta^{1+ {\delta}\mu_{L_\delta}}},
$$
which is exactly \eqref{e2.4}. \hfill$\Box$

\begin{rem} Using $f(v)\leq v$, and the solution $\overline v(t,x)$ of \eqref{e2.1} with  $f(v)=v-g(v)$ replaced by $v$, we obtain the (sharper) converse inequality
$$\frac{e^{T_\e}}{T_\e^{3/2}}\geq\frac{C}\e, 
$$
thus an asymptotic expansion of $T_\e$:
$$
T_\e={\mathrm {ln}}\frac1\e+\frac32{\mathrm {ln}}{\mathrm {ln}}\left(\frac1\e\right)+o_{\e\to0}\biggl({\mathrm {ln}}{\mathrm {ln}}\left(\frac1\e\right)\biggl).
$$
\end{rem}
 
\section{Communications between the road and the field}
The goal of this section is to prove that, if the solution on the road is of a certain order at some time and on a certain interval, then the solution in the field will
be of the same order, possibly in a square with a smaller size and a little later in time. We also want to prove that the converse holds: if the solution is of some order at some time and some point in the field, this is transmitted to the road. Such results can be seen as weak versions of the Harnack inequality (a bound at a certain time and point would entail the same bound in a whole neighbourhood, possibly at later times) but this will be sufficient for our purpose. See \cite{BDR}
 for estimates that are more in the spirit of the Harnack inequality.
\begin{lm}
\label{l2.1}
Consider $t_0\geq 1$, $x_0\in\RR$, $L\geq1$ and $\e>0$ (not necessarily small) such that 
$$
u(t_0,x)\geq\e\ \hbox{on $[x_0-L,x_0+L]$.}
$$
There is $c_L>0$ (universal otherwise) such that 
$$
v(t,x {,y})\geq c_L\e\ \hbox{on $[t_0+1,t_0+2]\times[ {x_0}-L, {x_0}+L]$ {$\times [0,1]$.}}
$$
\end{lm}
\noindent{\bf Proof.} Without loss of generality, we may translate time and space so as to have
$
t_0=1,\ x_0=0.$ Notice then that, because $v(t,x,0)\geq0$ we have
$$
u_t+(-\partial_{xx})^\alpha u+\mu u\geq0,\quad t\geq1,\ x\in\RR.
$$
Recall that the fundamental solution of the fractional heat equation of order $\alpha$, that we call $G_\alpha(t,x)$, is uniformly bounded away from 0 on $[1,2]\times [-L-1,L+1]$. Thus 
$$
u(t,x)\geq e^{-\mu  {(t-1)}}\int_{\vert y\vert\leq1}G(t,x- {x'})u(1, {x'})d {x'},\ \hbox{for $t>1$ and $\vert x\vert\leq L$.}
$$
This implies 
$$
u(t,x)\geq c_L\e,\ \hbox{for $t\in[1,2]$ and $\vert x\vert\leq L$.}
$$ 
Then, recall that $\di\frac{f(v)}v$ is bounded from below - say, by $-\Lambda>0$, and that
$$\tilde v(t,x,y)=e^{\Lambda t}v(t,x,y)
$$ is a super-solution to the heat equation, while the boundary condition reads
$$
\partial_y\tilde v+\nu \tilde v\geq c_L\e,
$$
for a possibly different $c_L$. Thus we have $\tilde v\geq\underline v$, where
$$
\begin{array}{rll}
(\partial_t-\Delta)\underline v=&0\quad (t\in(1,2],x\in\R,y>0)\\
\partial_y\underline v+\nu \underline v=&c_L\e\un_{[-2L,2L]}(x),\quad{y=0}\\
\underline v({1},x,y)=&0.
\end{array}
$$
By parabolic regularity we have, for some universal ${C>0}$:
$$
\vert\nabla \underline v(t,x,y)\vert\leq C\e,\quad t\in[1,2],\vert x\vert\leq 3L/2,0\leq y\leq1.
$$
Thus, there is $y_0\in(0,1)$, independent of $x_0$ - and thus of $\e$ such that
$$
{\underline v}(t,x,y_0)\geq \frac{c_L}2\e,\quad 1/2\leq t\leq 1,-x_0-3L/2\leq x\leq x_0+3L/2.
$$
And the classical parabolic Harnack inequality implies the lemma. Note that, due to \cite{BCN}, Section 3, one may push it to the boundary thanks to the Robin condition, at the expense of considering
$
\tilde{\underline v}(t,x,y):=e^{\nu y}\underline v(t,x,y). 
$\hfill$\Box$

\begin{lm}
\label{l2.2}
Consider $t_0\geq 1$, $x_0\in\RR$, and $\e>0$ such that 
$$
v(t_0,x_0,1)\geq\e.
$$
For all $L>0,$ there is $c_L>0$ (universal otherwise) such that 
$$
u(t,x),v(t,x,y)\geq c_L\e\ \hbox{on $[t_0+1,t_0+2]\times[x_0-2L,x_0+2L]\times[{0,1} ]$.}
$$
\end{lm}
\noindent{\bf Proof.} Once again there is no loss in generality by assuming $t_0=1$, $x_0=0$. The classical Harnack inequality applied to $v$ entails a lower bound of the order $\e$ at least for $v(t,x,1)$ for $t\in[1,2]$ and $-2L\leq x\leq 2L$. Fix now $L>0$,  for all $\delta\in(0,1)$ there is $c_\delta>0$ (we omit the dependence in $L$) such that 
$$
v(t,x,y)\geq c_\delta\e,\quad (t,x,y)\in[1,2]\times[-L,L]\times[ {\delta,1} ].
$$
Assume the existence of $x_1\in[-L,L]$ and $t_1\in[1,2]$ such that $v(t_1,x_1,0)$ is much smaller than its order of magnitude in the field. This is equivalent to assuming the existence of a sequence of solutions $(u_n,v_n)$ of \eqref{e1.1}, such that the following situation holds:
\begin{itemize}
\item for $t\in[1,2]$, $x\in[x_1-L/2,x_1+L/2]$ and $y=-1$, then $v_n(t,x,y)\geq c\e$ (dependence on $L$ omitted),
\item there is $t_1\in[1,2]$ such that $v_n(t_1,x_1,0)\leq1/n$.
\end{itemize}
Remember that $v_n$ is uniformly bounded from above. So, by parabolic regularity, (a subsequence of) the sequence $(u_n,v_n)_n$ converges, on $[1,2]\times[x_1-L/2,x_1+L/2]\times [-1,0]$ to a limiting function $(u_\infty,v_\infty)$ which is not identically equal to 0 due to the first assumption on $v_n$. The Hopf Lemma implies
$\partial_yv_\infty(t,x,0)<0$, thus the exchange condition yields
$$
\mu u_\infty(t_1,x_1)-\nu v_\infty(t_1,x_1,0)<0.
$$
This contradicts the fact that $v_\infty(t_1,x_1,0)=0$. Now, we have $u(t,x)\geq\underline u(t,x)$, with
$$
\left\{
\begin{array}{rll}
\underline u_t+(-\partial_{xx})^\alpha \underline u+\mu \underline u=&  {\nu} c\e\un_{[-L,L]}(x), \quad t>0,\ x\in \R\\
\underline u(0,x)=&0.
\end{array}
\right.
$$
Thus, for $t\in[1,2]$ we have
$$
u(t,x)\geq c\e {e^{-\mu t}}\int_0^1\int_{-L}^LG_\alpha(t-s,x- {x'})d {x'}ds\geq c'\e
$$
for a constant $c'$ that only depends on $L$. \hfill$\Box$

\section{Bounds to the full model}    
The starting point of the analysis is the (computationally non trivial) result, whose main line of the proof are
given in \cite{BCRR}, and proved in full length in \cite{ACC}.
\begin{thm}
[\cite{ACC}, Chapter 4]
\label{t3.1} Let $(\overline u(t,x),\overline v(t,x,y))$ solve
\begin{equation}
\label{e3.1}
\left\{
\begin{array}{rcll}
 \partial_{t}\o v-\Delta\o v &=&\o v,&\quad t>0, x \in\R, y\red{>}0\\
 \partial_{t}\o u+(-\partial _{xx})^{\alpha}\o u&=& -\mu \o u +\o v- k \o u,&\quad t>0, x \in \R\\
 \partial_{y}\o v&=&\mu\o  u - v, & x \in t>0, \R,y=0,\\
\end{array}
\right.
\end{equation}
with $(\o u(0,x),\o v(0,x,y))=(u_0(x),0)$ and 
$u_0\not\equiv 0$ nonnegative and compactly supported. There exists a function $R(t,x)$
and constants $\delta>0$, $C>0$ such that
\begin{enumerate}
\item we have, for large $x$:
\begin{equation}
\label{e3.2}
\biggl\vert \o u(t,x)-\frac{8\alpha\mu \sin(\alpha\pi)\Gamma(2\alpha)\Gamma(3/2)}{\pi}\frac{e^{t}}{t^{3/2}\vert x\vert^{1+2\alpha}}\biggl\vert\leq R(t,x),
\end{equation}
\item and the function $R(t,x)$ is estimated as
$$
0\leq R(t,x)\leq C\left(e^{-\delta t}+\frac{e^{t}}{\vert x\vert^{\min(1+4\alpha,3)}}+\frac{e^{t}}{\vert x\vert^{1+2\alpha}t^{\frac{5}{2}}}\right).
$$
\end{enumerate}
\end{thm}
Note that this result readily entails the upper bound in Theorem \ref{t1.1}, so that it suffices to prove the lower bound. We first prove it when $g$ is compactly supported in $(0,1]$, then indicate the necessary changes for a general $g$. \subsection{Proof of Theorem \ref{t1.1} when $g$ is compactly supported in $(0,1]$}
Let us pick $\lambda\in({\e},1/\mu)$ and $x_0>0$ (the same argument would apply for $x_0<0$) very large, we set
$$
u(1,x_0):=\e.
$$
From Theorem \ref{t3.1}, applied at time $t=1$, the function $u(1,x)$ is of the order $\e$ (and also of the order $1/x_0^{1+2\alpha}$)
on any interval around $x_0$ whose length will not exceed, say, $\sqrt{x_0}$. Thus, from Lemma \ref{l2.1}, we have
\begin{equation}
\label{e3.50}
v(1,x,y)\geq c\e\ \hbox{for $(x,y)\in[x_0-\sqrt x_0,x_0+\sqrt x_0] {\times [0,1]}$.}
\end{equation}

\noindent We ask how much time it will take for $u$ to reach the value $\lambda$ at $x_0$.
From \eqref{e3.2} we have, as soon as $\e{<\theta}$ is small enough - that is, if $x_0$ is large enough - and for all $L>0$:
$$
u(1/2,x)\leq c_L\e\quad \hbox{for all $x\in[-x_0-2L,x_0+2L]$}.
$$
Then, translate the point $(x_0,0)$ to the origin, and let this time 
$\underline v(t,x,y)$ solve the two-dimensional Fisher-KPP equation with Dirichlet conditions on the road:
\begin{equation}
\label{e3.4}
\left\{\begin{array}{rll}
(\partial_t-\Delta-1)v+g(v)=&0, \quad t\geq 1,x\in\R,y{>}0\\
\underline v(t,x,0)=&0,  \quad {t\geq 1, \ x\in\R}\\
\underline v(1,x,y)=& {c}\e\un_{[-\sqrt x_0,\sqrt x_0]}(x)\un_{[0,1]}(y).
\end{array}\right.
\end{equation}
As long as $\underline v\leq\theta$ everywhere, it solves the linear equation, that is, with $g(v)=0$. 
In such a case it consists of the product of two solutions of the heat equation times the exponential:
\begin{equation}
\label{e3.6}
\underline v(t,x,y)=\frac{\underline v^{1D}(t,y)}{\sqrt\pi}\biggl(\int_{ {[}-\sqrt x_0,\sqrt x_0{]}}e^{-(x-x')^2/4t}dx'\biggl),
\end{equation}
where $\underline v^{1D}(t,y)$ is the solution of the Dirichlet heat equation in $y$, and is exactly given by \eqref{e2.6}. The function $v^{1D}$ reaches $\theta$ in
a time  $T^1_\e$ given by equation \eqref{e2.10}, in other words
$$
T_\e^1=O\left({\mathrm {ln}}\left(\frac1\e\right)\right).
$$
This time is too short for the solution of the heat equation in $x$ to decay significantly on, say, the interval $[-L,L]$ with $L$ large but finite (any size $L$ which is an $o(\sqrt x_0)$ will do). We have indeed, for $\vert x\vert\leq L$ and $t\leq T_\e^1$:
$$
\begin{array}{rll}
\di\frac1{\sqrt t}\di\int_{-\sqrt x_0}^{\sqrt x_0}e^{-\frac{(x-x')^2}{4t}}dx'
\geq&  {C}\di\int_{-\frac{(\sqrt x_0+L){^2}}{2\sqrt t}}^{\frac{(\sqrt x_0-L) {^2}}{2\sqrt t}}e^{-\xi^2}d\xi\\
\sim& C,
\end{array}
$$
simply because $x_0\sim\e^{-1/(1+2\alpha)}$ and $T^1_\e$ is of the order $\di{\mathrm{ln}} {\left(\frac1\e\right)}$. Thus,  the function $y\mapsto \underline v(T_\e^1,0,y)$ reaches maximum of the order $\theta$, at a point of the order $\sqrt{T_\e^1}$, while it is of the order $\di\frac1{\sqrt{T_\e^1}}$ for $y\sim1$. Then, we run \eqref{e3.4} again, from $T_\e^1$, with 
$$
\underline v(T_\e^1,x,y)=\un_{[{-}L',L']}(x)\un_{[ {1,L'} ]}\frac{e_{L'}(x,y)}{\sqrt {T_\e^1}},
$$
the function $e_{L'}(x,y)$ being the first eigenfunction of the Dirichlet Laplacian in the rectangle $(-L',L')\times( {1,L'} )$, and $L'$ a large number (in fact we could take $L'=L$). We have
$$
e_{L'}(x,y)=\sin\biggl(\frac{\pi}{2L'}x\biggl)\sin\biggl(\frac{\pi}{L'-1}(y-1)\biggl),
$$
the first eigenvalue being still an $O\left(\di\frac1{{L'}^2}\right)$. And so, for a time $T_\e$ given by \eqref{e2.12}, there is $q>0$ independent of $\e$ such that $\underline v(T_\e,0,1)\geq q$.

To conclude the proof, it remains to prove the existence of $q'>0$ universal such that $u(T_\e,0)\geq q'$. This is, however, easy: because $x_0$ is arbitrary, we have
$$
v(T_\e,x,1)\geq q\quad\hbox{for $x_0-L\leq x\leq x_0+L$,}
$$
and Lemma \ref{l2.2} implies the desired bound for $u$. 
\subsection{The lower bound for a general concave nonlinearity}    
We write again
$$
f(v)= v-g(v),\quad g(0)=g'(0)=0,\ g(1)=1,\ g''>0 \text{ on [0,1]}.
$$
Thus we have, for all $v\in[0,1]$: $g(v)=O(v^2)$, and this is what we will really use. In view of what we have already done when $g$ vanishes in a vicinity of 0,
what we really have to do is study the function $v(t,x,y)$ solving
\begin{equation}
\label{e4.1}
\begin{array}{rll}
v_t-\Delta v- v=&g(v)\quad(t>0,x>0,y {>}0)\\
v(t,x,0)=&0\\
v(0,x,y)=&c\e\un_{[-\sqrt x_0,\sqrt x_0]}(x)\un[-1,0](y),  
\end{array}
\end{equation}
with $\e=\di\frac1{1+x_0^{1+2\alpha}}$. In view of the proof of Theorem 3.1, and Section 4, the main property that we have to prove is the following.
\begin{lm}
\label{l4.1}
Let $T_\e^1$ be given by \eqref{e2.10}. There is $q>0$ universal such that
\begin{equation}
\label{e4.2}
v(T_\e^1,0,1)\geq \frac{q}{\sqrt {T_1^\e}}. 
\end{equation}
\end{lm}
It then suffices, as in the preceding section, to put $v$ above the solution $\underline v_{\delta,L,L'}$ of 
$$
\left\{\begin{array}{rll}
\biggl(\partial_t-\Delta-(1-\delta)\biggl)v_{\delta,L,L'}=&0\quad (t>T_\e^1,-L<x<L, {1}<y< {L'})\\
v_{\delta,L,L'}(0,x,y)=&\di\frac{q}{\sqrt {T_\e^1}}e_{L,L'}(x,y),
\end{array}\right.
$$
with $\delta$ small, and $e_{L,L'}(y)$ the first eigenfunction of the Dirichlet Laplacian in $(-L,L)\times {(1,L')} $. At time 
$$
T_\e^2=T_\e^1+\biggl(\frac12+O(\delta)+O\left(\frac1{L^2}\right)+O\left(\frac1{{L'}^2}\right)\biggl){\mathrm{ln}}T_\e^1,
$$
we have $\underline v_{\delta,L,L'}(t,x,y)\geq C\delta$ on $(-L,L)\times {(1,L')} $, and one finishes the proof of Theorem \ref{t1.1} by Lemmas \ref{l2.1} and \ref{l2.2}. 

Let us therefore present the

\noindent {\bf Proof of Lemma \ref{l4.1}.} Call $X=(x,y)\in\R\times\R_{+}$ a generic point of the upper half-plane $\R\times\R_-$. Let $G(t,X)$ be the fundamental solution of the Dirichlet heat equation in the upper  half-plane, we have
$$
G(t,X,X')=G_0(t,y,y')G_1(t,x),
$$
the function $G_0$ being the Dirichlet fundamental solution (see \eqref{e2.6}) whereas $G_1$ is the standard Gaussian 
$G_1(t,x)=\di\frac{e^{-x^2/4t}}{\sqrt{4\pi t}}$. 
The Duhamel formula
yields
$$
v(t,X)=\e e^t\int_{\R^2_{+}}G(t,X,X')v(0,X')dX'-\int_0^t\int_{\R^2_{+}}e^{t-s}G(t-s,X,X')g(v(s,X'))dsdX'.
$$
We call $e^{t-s}D(t-s,s,X, X')$ the integrand of the second integral in the right handside of the above inequality; we have
$$
v(s,X)\leq  {e^s\e}\int_{\R^2_{+}}G(s,X,X')v(0,X')dX'\leq C\e{^2e^s}\int_{-1}^0G_0({s},y,y')dy'.
$$
Because $g(v)=O(v^2)$ we get, taking \eqref{e3.6} and \eqref{e2.5} into account: 
$$
D(t-s,s,X,X') \leq C\e^{4}G(t-s,X,X')\frac{{y'}^2e^{-2(y-y')^2/5s}}{s(1+s^2)} {e^{2s}},
$$
$C>0$ universal. Note that we have only estimated the integral for $s\geq1$, the integral for $s\leq 1$ being negligible. Integrating in $x'$ and specialising at $x=0$ we get
$$
\int_0^t\int_{\R^2_{+}}e^{t-s}G(t-s,X,X')g(v(s,X'))dsdX'\leq C\e^{ {4}}\int_0^t\int_{\R_{+}}e^{t {+}s}E(t-s,s,y')dsdy',
$$
where
\begin{equation}
\label{e4.4}
E(t-s,s,y')\lesssim \frac{\vert 1+ y'\vert^3 e^{-(1+y')^2/4(t-s)}e^{-{y'}^2/3s}}{(1+\vert t-s\vert)(t-s)^{1/2}s(1+s^2)}.
\end{equation}
We are going to prove the inequality 
\begin{equation}
\label{e4.3}
\int_0^{T_1^\e} {\int_{\R\times \R_+}}e^{T_1^\e-s}D(T_1^\e-s,s,(0,1),X')dsdX'\leq C(\sqrt\e+\frac{1}{{\mathrm{ln}}\frac1\e})\frac1{\sqrt{T_1^\e}},
\end{equation}
$C>0$ universal. There is nothing special about the point $X=(0,1)$ the inequality would be valid for all neighbouring points, at the expense of increasing $C$. 
Recall the inequality (see Sections 2 and 3) for $v$:
$$
v^{1D}(T_1^\e,1)\geq\frac{q}{\sqrt T_1^\e}.
$$
This, combined to \eqref{e4.3}, will imply the lemma.
As we will set, eventually, $t=T_1^\e$,  we will always assume
$$
t=O({\mathrm{ln}}\frac1\e).
$$
we cut the time interval $(0,t)$ into two.

\noindent{\bf 1. $s\in(0,\kappa t)$, $\kappa {>0}$ small.} We will use the factor $e^{-{y'}^2/3s}$ to make the integral convergent, and make the change of variables $y'\mapsto z'=y'/\sqrt s$.Thus we have
$$
E(t-s,s,y')\lesssim\di\frac{1+\vert z'\vert^3e^{-{z'}^2/3}}{(1+(t-s))\sqrt{t-s} s},
$$
using $e^{t+s}\leq e^{(1+\kappa) t}$ and the definition of $T_1^\e$ we end up with
$$
\e^2\int_0^{\kappa t} {\int_{\R\times \R_+}}e^{t-s}D(t-s,s,(0,1),X')dsdX'\lesssim \frac{\e^{4}te^{(1+\kappa)t}}{t^{3/2}}.
$$
And so,
\begin{equation}
\label{e4.5}
\e^2\int_0^{(1+\kappa){T_1^\e}}{\int_{\R\times \R_+}}e^{T_1^\e}D(T_1^\e,s,(0,-1),X')dsdX'\lesssim\frac{\e^{1-\kappa}}{\sqrt {T_\e^1}}.
\end{equation}
\noindent{\bf 2. The range $s\geq\kappa t$.} This time we rely on the part $e^{-\frac{(1+y')^2}{4(t-s)}}$ to make the spatial integral convergent, and we will have to be a little careful about the $e^{t+s}$ factor. As for the powers ${y'}^2$, we dominate them by $1+\vert 1+y'\vert^2$. We make the change of variables $y'\mapsto\di\frac{z'}{\sqrt{t-s}}$ and we have
$$
\begin{array}{rll}
e^{t-s}D(t,s,(0,-1),X')\lesssim&\e^2G_1(x')\frac{e^{t+s}e^{-\frac{{z'}^2}{4}}\sqrt{t-s}}{(1+\vert z'\vert^3)s^3}\\
\lesssim&\e^2G_1(x')\frac{e^{t+s}e^{-\frac{{z'}^2}{4}}}{t^{5/2}}.
\end{array}
$$
Integrating on $(0,t)\times\RR\times\RR_+$ yields
$$
\e^2\int_0^t\int_{\RR\times\RR_+}e^{t-s}D(t,s,(0,-1),X')dsdX'\leq \frac{Ce^{2t}}{(T_1^\e)^{3/2}}\leq\frac{C}{{\mathrm{ln}}\frac1\e}.\frac1{\sqrt {T_1^\e}}.
$$
Putting everything together yields \eqref{e4.3}, hence the lemma. \hfill$\Box$

\noindent {\bf Acknowledgement.} The research leading to these results has received funding from the ERC under the European Union's Seventh Frame work Programme (FP/2007-2013) / ERC Grant Agreement 321186 - ReaDi.

\noindent 
{\footnotesize
 
}
\end{document}